\documentclass[11pt]{amsart}

\usepackage{amsmath,amsthm,amssymb}

\newtheorem{theorem}{Theorem}
\newtheorem{definition}[theorem]{Definition}
\newtheorem{proposition}[theorem]{Proposition}
\newtheorem{corollary}[theorem]{Corollary}

\newcommand{\N}{\omega^\omega}
\newcommand{\dom}{\operatorname {dom}}

\title{A note on the diamond operator}

\author[L.~Westrick]{Linda Westrick} 
\address{Penn State University Department of Mathematics, University Park, PA}
 \email{westrick@psu.edu}

\thanks{This work was done while the author was in attendance at Dagstuhl Seminar 18361, ``Measuring the Complexity of Computational Content: From Combinatorial Problems to Analysis.''  The author was also supported by the Cada R. and Susan Wynn Grove Early Career Professorship in Mathematics, and partially supported by NSF grant DMS-1854107.  We thank Arno Pauly for many useful and clarifying discussions. We also thank the anonymous referees, who provided additional context including Corollary 2, and uncovered some minor bugs.}

\begin{document}

\begin{abstract}
We show that if $1 \leq_W F$ and $F \star F \leq_W F$, then $F^\diamond \leq_W F$, 
where $\star$ and $\diamond$ are the following operations 
in the Weihrauch lattice: $\star$ is the compositional product, which allows 
the use of two principles in sequence, while the diamond operator
$\diamond$ allows an arbitrary but finite number of uses of the given principle in sequence.  This 
answers a question of Pauly.
\end{abstract}

\maketitle

Informally, let $F$ be a problem such as 
``given a continuous function on $[0,1]$, find an input where it obtains its maximum value''.  Although this problem 
is well-posed, in general there is no constructive method for solving it.  However, if we add
a single-use
oracle access to a solver for the problem $G$: ``given an infinite binary tree, find a path through it'',
$F$ becomes constructively solvable.  Formally, one may define the notion of Weihrauch 
reducibility to express this relationship: we write $F \leq_W G$ and say that $F$ is 
Weihrauch reducible to $G$.  In what follows we assume familiarity with the Weihrauch 
lattice; for definitions and references, see the recent survey \cite{BrattkaGherardiPauly}.

Sometimes a problem $F$ becomes solvable if we allow single-use oracle access to solvers 
for multiple problems 
$G$ and $H$, or if we allow the oracle multiple uses of $G$.  In this case we distinguish 
between the different ways that the oracle can be used.  The compositional product,
roughly speaking, takes two problems $G$ and $H$ and returns the problem 
$H \star G$: ``use $G$ and then $H$''.  The diamond operator, roughly speaking,
takes a problem $F$ and returns the problem $F^\diamond$: ``use $F$ an 
arbitrary but finite number of times''.  The formal definitions are given below.  A natural 
question was asked by Pauly \cite{Dagstuhl2018}: 
suppose that given a single-use oracle access to $F$, we can solve $F\star F$. 
Does it follow that a single-use oracle access to $F$ is enough to solve $F^\diamond$?
We show that it is.  

\begin{theorem}\label{thm1}
  Suppose $F$ is a Weihrauch problem with $1\leq_W F$ and $F \star F \leq_W F$.
  Then $F^\diamond \leq_W F$.
  \end{theorem}

It follows that for any $F$ with $1 \leq_W F$, $F^\diamond$ is the least Weihrauch 
degree above $F$ that is closed under $\star$.

\begin{corollary}
If $F$ and $G$ are Weihrauch problems with $1 \leq_W F \leq_W G$, 
and if $G\star G \leq_W G$, then $F^\diamond \leq_W G$.
\end{corollary}

The technical condition $1\leq_W F$ in Theorem \ref{thm1}
removes a trivial counterexample by 
ensuring that $F$ has a computable instance.  If $F$ does not have a
computable instance, then we cannot have $F^\diamond \leq_W F$ 
because $F^\diamond$ has computable instances (the input can request to use $F$ 
zero times), while a formal requirement of the $\leq_W$ definition is that 
the oracle is used exactly once, requiring an input to $F$ to be produced.  

The diamond operator was first introduced by Neumann and Pauly 
\cite{NeumannPauly2018} for the
purpose of comparing the complexity of the Blum-Shub-Smale (BSS) 
computation model and the Type-2-Effectivity (TTE) model.  
In the BSS model, a computation has access to an arbitrary
finite number of uses of the relations $<$ and $=$ on real numbers, 
which is equivalent to access to $LPO^\diamond$.

As Neumann and Pauly note in their paper, the diamond operator 
was also essentially defined by Hirschfeldt and Jockusch 
\cite{HirschfeldtJockusch2016} using a game characterization.  
They define a relation $F \leq_\omega G$ to mean that whenever 
$G$ holds in an $\omega$-model of second order arithmetic, 
so does $F$.  Given any pair of principles $(F,G)$, they define 
a two-player game such that Player II has a winning strategy 
if and only if $F \leq_\omega G$.  Briefly, the game is as follows: 
Player I plays an instance $X_0$ of $F$; Player II may either play
a solution $S_0 \leq_T X_0$ (in which case she wins), 
or an instance $Y_0$ of $G$
such that $Y_0 \leq_T X_0$.  Then Player I must play a solution 
$T_0$ to $Y_0$.  Player II may now either play a solution
$S_0$ to $X_0$ with $S_0 \leq_T X_0\oplus T_0$ (and win), 
or a new instance $Y_1$ of $G$ with $Y_1 \leq_T X_0 \oplus T_0$.
The game continues in a similar fashion, with Player I winning 
if Player II never produces a solution to $X_0$.
Although $\leq_\omega$ is characterized by this game, 
the strategies are not required to be effective.  So Hirschfeldt 
and Jockusch also define a relation
$F\leq_{gW} G$, or  ``$F$ is generalized Weihrauch reducible to $G$''
to mean that  $F \leq_\omega G$ and furthermore Player II has a 
\emph{computable} winning strategy in the
game.  As noted by Pauly \cite{Dagstuhl2018}, it is a matter of 
definition-chasing to see that
for all problems $F$ and $G$, we have
$F \leq_{gW} G$ if and only if $F \leq_W G^\diamond$.

Now we give the definitions. Let $\Phi$ be a universal functional.

\begin{definition}\label{def2}
  If $F,G:\subseteq \omega^\omega \rightrightarrows \omega^\omega$
   are Weihrauch problems then the \emph{compositional product} 
   $F\star G$ is the problem whose
  domain is $$\{(x,e) \in \N\times \N : x \in \dom G \text{ and for all }
  y \in G(x), \Phi(e,y) \in \dom F\}$$
  and whose output is a pair $(y,z)$ such that $y \in G(x)$ and
  $z \in F(\Phi(e,y))$.
  \end{definition} 
  
  The compositional product was first defined in \cite{BrattkaPauly2018}, 
  and more details about it can be found there. The most 
  accessible definition is found in \cite{BrattkaGherardiPauly},
  where $(F \star G)(x,p)$ is defined as 
  $\langle id \times F \rangle \circ \Phi_p \circ G(x)$.  In that definition 
  $\langle id \times F\rangle$ is the problem which, on input 
  $(y,w)$, outputs $(y,F(w))$; thus, this problem interprets the output 
  of $\Phi_p$ as a pair.  
  
 Our definition may appear to differ slightly from the one given in 
 \cite{BrattkaGherardiPauly}.
  To see that the two definitions lie in 
  the same Weihrauch degree, first see that if we 
  are given $(x,p)$, we may let $e$ be the functional 
 defined by letting $\Phi(e,y)$ be the second component of 
 $\Phi(p,y)$.  Then given $(y,z)$ with
  $y \in G(x)$
   and $z \in F(\Phi(e,y))$, we may compute $w$ 
   as the first component of $\Phi(p,y)$, and then 
   check that $(w, z) \in 
   \langle id \times F \rangle \circ \Phi_p \circ G(x)$.
   On the other hand, if we are given input $(x,e)$, 
   we may define $p$ to be the functional which 
   on input $y$, outputs $(y,\Phi(e,y))$.  Then if 
   $(y,z) \in \langle id \times F \rangle \circ \Phi_p \circ G(x)$,
   we have $y \in G(x)$ and $z \in F(\Phi(e,y))$ as needed.

Now we define the diamond operator.
Let $\Phi^\ast$ be a universal functional with a space for a
Weihrauch problem plug-in: 
given a Weihrauch problem $F$ and any $d \in \omega^\omega$ as input,
in addition to its usual algorithmic operations involving $d$ the functional may
at any time
ask for an answer to $F(\Phi(t))$, where $t$ is the contents
of all the tapes at the time of the question.  A new tape is
created which contains an answer (assuming that $\Phi(t)$ is
total and in the domain of $F$; failure of either of these
conditions results in an infinite loop).  If the functional
halts, its output is considered to be $\Phi(t)$,
where $t$ is the contents of all the tapes.  If $\Phi(t)$ 
fails to be total, then the computation is said to have no output.
We let $\Phi^\ast(F,d)$ denote a run of this process 
using Weihrauch problem $F$ and starting with input
$d\in \N$.  Here we think of $d$ as describing both what 
$\Phi^\ast$-program to run and any input data.   The 
domain of $F^\diamond$ is the set of all $d \in \omega^\omega$
for which this functional always halts, no matter what answers 
to the $F$-questions are supplied.

\begin{definition}
  If $F$ is a Weihrauch problem then $F^\diamond$ is the problem whose domain
  is $$\{d \in \N: \Phi^\ast(F,d) \text{ always halts with total output}\}$$
  and whose output is the output of $\Phi^\ast(F,d)$.
\end{definition}

In the definition of $\star$, $\Phi$ was specified 
as a universal functional, but we did not get into any technical details about 
how $\Phi$ expects its inputs, or what it would do if given 
multiple tapes as input.  Let us be more specific because it also 
allows us to specify a technical detail of $\Phi$'s implementation 
that will come in very handy later.
We declare that $\Phi$'s single-tape mode of operation is as follows:
it expects an input of the form $(\dots(((n,x_1),x_2),\dots),x_r)$,
where $n \in \mathbb N$ and there are finitely many $x_i \in \N$.  
It parses the 
data, finds the number $n$, interprets 
this number as a program description, and then applies that 
program to the rest of the data $x_1,\dots,x_r$.  If we apply $\Phi$ 
in a multi-tape situation, with finitely many tapes occupied with 
data $d_1,\dots, d_r$, we declare the result to be the same as applying the 
single-tape version of $\Phi$ to the input $(\dots((d_1(0),d_1),d_2),\dots,d_r)$.
The reason for wanting to distinguish a particular number as the 
program number will be made clear below as we will designate 
that number using the recursion theorem.  Going forward, we write 
$(n,x_1,\dots,x_r)$ instead of the more cumbersome
$((\dots((n,x_1),x_2),\dots,x_r))$.\footnote{Yes, $\Phi$ diverges 
if given an infinite regress of pairings as input, but this situation will not arise.}

When comparing $\star$ and $\diamond$, it is tempting to hope that 
any problem reducible to $F^\diamond$ is actually reducible to 
$(F\star(F\star \dots  (F\star F)\dots))$ for some finite number of applications
of $\star$.  If this were true, Theorem \ref{thm1} would immediately follow.  
Failing that, one might hope that at least for any fixed instance $d$ of 
$F^\diamond$, there might be a bound on the number of calls to $F$ 
in $\Phi^\ast(F,d)$.
However,
the following example due to Pauly and Yoshimura \cite{Dagstuhl2018}
shows that this need not be.  Let $(q_i)_{i\in \omega}$ be 
a sequence of elements of $\N$ such that for all $i$, $q_i \not\leq_T \bigoplus_{j\neq i} q_j$.
Define $F$ by $F(0^\omega) = \{iq_i : i \in \omega\}$ and $F(q_{i+1}) = q_i$.
Let $q_0$ be the problem of producing $q_0$ given any element of $\omega^\omega$.
Then although $q_0 \leq_W F^\diamond$, there is no bound on the number of 
uses of $F$ that may be required, even for the fixed input $0^\omega$.

A key observation in the proof of Theorem \ref{thm1} is that the recursion theorem 
relativizes.  That is, the usual proof of the recursion theorem also proves this:
\begin{proposition}[Relativized recursion theorem]
For any computable function $f$, there is a number $n$ such that for all 
oracles $X$, $\Phi^X_{f(n)} \cong \Phi^X_n$.
\end{proposition}

Finally we prove the theorem.
\begin{proof}[Proof of Theorem 1.]  Since $F\star F \leq_W F$,
  let $\Delta$ and $\Gamma$ be functionals such that for all
  $(x,e) \in \dom(F\star F)$, $\Delta(x,e) \in \dom F$ and for all
  $w \in F(\Delta(x,e))$, $\Gamma(w,(x,e)) = (y,z)$, where $y \in F(x)$
  and $z \in F(\Phi(e,y))$.

  Recall our assumption that $\Phi$ looks at the first bit of
  the innermost real of
   its input to
  determine what program to run.
  By the recursion theorem, let $n \in \omega$ be a fixed point such that
  for all oracles $(d,y_1,\dots,y_k) \in (\N)^{<\omega}$, the computation
  $$\Phi(n,d,y_1,\dots,y_k)$$
  does the following:
  \begin{enumerate}
    \item Starts simulating $\Phi^\ast(\cdot, d)$ without
      a Weihrauch problem plugged in.
    \item Upon the $i$th oracle request,
      (up to $k$ requests), $y_i$ is made available
      to the simulation.
    \item\label{l3} If the
      simulation halts after $k$ or fewer requests, output a fixed
      computable element of $\dom F$.
    \item\label{l4} If the simulation makes a
      $(k+1)$st request, halt and return
      $$\Delta(\Phi(t_{k+1}),(n,d,y_1,\dots,y_k)),$$
      where $t_{k+1}$ is the contents of the simulated 
      tapes at the time of the $(k+1)$st
      request.
  \end{enumerate}
  To be explicit about computational inputs and outputs,
  let us give the name $\Psi$ to a fixed functional such that 
  if the above computation gets to stage 4 on input 
  $(d,y_1,\dots,y_k)$, then 
  $\Psi(d,y_1,\dots,y_k) = \Phi(t_{k+1})$ above.  If 
  the computation does not get to stage 4,  
$\Psi(d,y_1,\dots,y_k)$ diverges.

  Now we describe how to reduce $F^\diamond$ to $F$.
  Given $d \in \dom F^\diamond$,
  we ask $F$ for an answer to $\Phi(n,d)$. To guarantee that this
  question is in the domain of $F$,
  we verify a superficially stronger (but also necessary) property.

  {\bf Claim.}  Suppose that $y_1,\dots y_k$ are such that
  $y_i \in F(\Psi(d,y_1,\dots y_{i-1}))$ for all $i\leq k$.  Then
  $\Phi(n,d,y_1,\dots,y_k)$ is total and in the domain of $F$.
  
  Consider the set of all $(y_1,\dots,y_k)$ such that the
  hypotheses of the claim are satisfied.  This set is a tree
  when ordered by extension
  (with possibly continuum-sized branching).
  Its members represent all possible ways that $F$ could
  supply its information to the computation
  $\Phi^\ast(F,d)$.  Since $d \in \dom F^\diamond$,
  it is guaranteed that no matter how $F$ supplies its
  answers, the computation only asks finitely many questions, 
  and then halts.
  Therefore, this tree is well-founded.  We prove the claim
  by induction on the well-founded rank of the tree.
  If $(y_1,\dots,y_k)$ is a leaf, the computation terminates
  without asking any more questions to $F$.  In this case,
  $\Phi(n,d,y_1,\dots,y_k)$ halts in step (\ref{l3}) above
  and returns an element of $\dom F$.  If $(y_1,\dots, y_k)$
  is not a leaf, the question $\Psi(d,y_1,\dots,y_k)$ is
  asked, and $\Phi(n,d,y_1,\dots,y_k)$ halts in step (\ref{l4}).
  We claim the value it gives is in $\dom F$.  First, the hypothesis 
  on $(d,y_1,\dots, y_k)$ guarantees we are on a correct 
  computation path, so $\Psi(d,y_1,\dots,y_k) \in \dom F$,
  and
  for any $y_{k+1} \in F(\Psi(d,y_1,\dots,y_k))$, by induction
  we know that $\Phi(n,d,y_1,\dots,y_k,y_{k+1}) \in \dom F$.
  It follows that $(\Psi(d,y_1,\dots,y_{k}),(n,d,y_1,\dots,y_k)) \in \dom (F \star F)$,
  and so applying $\Delta$, we see $\Phi(n,d,y_1,\dots,y_k) \in \dom F$.
  This proves the claim.

  When $k=0$, the claim implies that $\Phi(n,d) \in \dom F$.
  Let $z_0 \in F(\Phi(n,d))$ be the answer given.  We now
  compute $F^\diamond(d)$ as follows.  Begin to
  simulate $\Phi^\ast(F,d)$.  If it ever halts, output what it
  outputs.   If it asks a first question $\Psi(d)$, then
  $\Phi(n,d)$ halted in case (\ref{l4}) and returned
  $\Delta(\Psi(d),(n,d))$.  By the Claim, $(\Psi(d),(n,d))$
  is a valid input to $F\star F$, so since
  $z_0\in F(\Delta(\Psi(d),(n,d)))$, we know that $\Gamma(z_0,(\Psi(d),(n,d)))$
  yields a pair $(y_1,z_1)$ such that $y_1\in F(\Psi(d))$
  and $z_1 \in F(\Phi(n,d,y_1))$.  Continue the simulation
  using $y_1$ as the simulated
  answer of $F$.  In general, we maintain that
  $y_k \in F(\Psi(d,y_1,\dots,y_{k-1}))$ and
  that $z_k \in F(\Phi(n,d,y_1,\dots,y_k))$.  Then if there is a
  $(k+1)$st question $\Psi(d,y_1,\dots,y_k)$, we run
  $\Gamma(z_k,(\Psi(n,d,y_1,\dots,y_k),(n,d,y_1,\dots,y_k)))$ to obtain
  $y_{k+1}$ which answers the question and $z_{k+1}$ which
  maintains the condition.  Eventually the simulation must stop 
  asking questions and halt with a correct simulated output.\end{proof}

\bibliographystyle{alpha}
\bibliography{bib}

\end{document}